\documentclass[12pt,a4paper,english]{article}
\usepackage{amssymb}
\usepackage{amsmath}
\usepackage{graphicx}
\usepackage[english]{babel}

\def\thebiblio#1{\subsection*{Referencias}\list
{[\arabic{enumi}]}{\settowidth\labelwidth{#1.}\leftmargin\labelwidth
 \advance\leftmargin\labelsep
 \usecounter{enumi}}
 \def\newblock{\hskip .11em plus .33em minus -.07em}
 \sloppy
 \sfcode`\.=1000\relax}

\input{tcilatex}

\begin{document}

\begin{center}
{\Large Forbidden patterns and shift systems}

\medskip

\noindent Jos\'{e} Mar\'{\i}a Amig\'{o}$^{\text{a}}$, Sergi Elizalde$^{\text{b}}$, Matthew B. Kennel$^{%
\text{c}}$
\end{center}

\noindent

\begin{center}
\noindent $^{\text{a}}$Centro de Investigaci\'{o}n Operativa,
Universidad Miguel Hern\'{a}ndez

03202 Elche, Spain


\noindent $^{\text{b}}$Department of Mathematics, Dartmouth College

Hanover, NH 03755-3551, USA

$^{\text{c}}$Institute for Nonlinear Science, University of
California, San Diego

La Jolla, CA 92093-0402, USA\

\bigskip
\end{center}

\bigskip


\noindent

\noindent \textbf{Abstract. }The scope of this paper is two-fold. First, to
present to the researchers in combinatorics an interesting implementation of
permutations avoiding generalized patterns in the framework of discrete-time
dynamical systems. Indeed, the orbits generated by piecewise monotone maps
on one-dimensional intervals have forbidden order patterns, i.e., order
patterns that do not occur in any orbit. The allowed patterns are then those
patterns avoiding the so-called forbidden root patterns and their shifted
patterns. The second scope is to study forbidden patterns in shift systems,
which are universal models in information theory, dynamical systems and
stochastic processes. Due to its simple structure, shift systems are
accessible to a more detailed analysis and, at the same time, exhibit all
important properties of low-dimensional chaotic dynamical systems (e.g.,
sensitivity to initial conditions, strong mixing and a dense set of periodic
points), allowing to export the results to other dynamical systems via
order-isomorphisms.

\medskip

Dynamical systems. Shift maps.


\bigskip

\section{Introduction}

Order has some interesting consequences in discrete-time dynamical systems.
Just as one can derive sequences of symbol patterns from such a dynamic via
coarse-graining of the phase space, so it is also straightforward to obtain
sequences of \textit{order patterns} if the phase space is linearly ordered.
It turns out that, under some mild mathematical assumptions, not all order
patterns can be materialized by the orbits of a given, one-dimensional
dynamic. Furthermore, if an order pattern of a given length is `forbidden',
i.e., cannot occur, its absence pervades all longer patterns in form of more
missing order patterns. This cascade of outgrowth forbidden patterns grows
super-exponentially (in fact, factorially) with the length, all its patterns
sharing a common structure. Of course, forbidden and allowed order patterns
can be viewed as permutations; allowed patterns are then those permutations
avoiding the so-called forbidden root patterns and their shifted patterns
(see Sect. 4 for an exact formulation). Let us mention at this point that
permutations avoiding generalized and consecutive patterns is a popular
topic in combinatorics (see, e.g., \cite{Babson,EliI,EliII}). It is in this
light that we approach order patterns in the present paper. In fact, the
measure-theoretical aspects of the underlying dynamical system play no role
in the combinatorial properties of the order patterns defined by its orbits
and hence will be only considered when necessary. Also for this reason we
will not dwell on the dynamical properties of shift systems and their role
as prototypes of chaotic maps once endowed with appropriate invariant
measures; see \cite{Baker,Doherty} for readable accounts.

Order relations belong rather to algebra than to continuous mathematics
because of their discrete nature. Only in the standard real line, order and
metric are coupled, leading to such interesting results as Sarkovskii's
theorem \cite{Sarko,Melo}. But even in this special though important
framework, order fails to be preserved by isomorphisms, that consistently
only address dynamical properties such as invariant measures, periodicity,
mixing properties, etc., and this reduces its applicability. Yet, order
relations have been successfully applied in discrete dynamical systems and
information theory, e.g., to evaluate the measure-theoretic and topological
entropies \cite{Bandt, Amigo1}. This paper is an extension of those
investigations. Isomorphisms that preserve the possibly existing order
relations of the dynamical systems they identify, are called
order-isomorphisms. The order isomorphy in one-dimensional dynamical systems
is the subject of \textit{Kneading Theory} \cite{Melo}. In this paper, we
will go beyond the framework of Kneading Theory in two respects: (i) the
maps need not be continuous (but piecewise continuous) and (ii) we will also
consider more general phase spaces (like finite-alphabet sequence spaces and
two dimensional intervals).

Forbidden order patterns, the only ones we will consider in this paper,
should not be mistaken for other sorts of forbidden patterns that may occur
in dynamics with constraints. Forbidden patterns in symbol sequences occur,
e.g., in Markov subshifts of finite type and, more generally, in random
walks on oriented graphs. On the contrary, the existence of forbidden
\textit{order} patterns does not entail necessarily any restriction on the
patterns of the corresponding symbolic dynamic: the variability of \textit{%
symbol} patterns is given by the statistical properties of the dynamic. As a
matter of fact, the symbolic dynamic of one-dimensional chaotic maps are
used to generate pseudo-random sequences, although all such maps used in
practice have forbidden order patterns. In general it is very difficult to
work out the specifics of the forbidden patterns of a given map, but we will
see that shifts on finite-symbol sequence spaces are an important exception:
the detailed analysis of the forbidden patterns of this transformations is
precisely the topic of this paper.

The existence of forbidden patterns is a hallmark of deterministic orbit
generation and thus it can be used to discriminate deterministic from random
time series. Indeed, thanks to the super-exponentially growing trail of
outgrowth forbidden patterns, the probability of a false forbidden pattern
in a truly stochastic process vanishes very fast with the pattern length
and, consequently, a time series with missing order patterns of moderate
length can be promoted to deterministic with virtually absolute confidence.
The quantitative details depend, of course, on the specificities of the
process (probability distribution, correlations, etc.). Only those chaotic
maps with all forbidden patterns of exceedingly long length seem to be
intractable from the practical point of view. Besides, applications need to
address some key issues, such as the robustness of the forbidden patterns
against observational noise, and the existence of false forbidden patterns
in \textit{finite}, random time series. We refer to \cite{Amigo3} for these
issues.

This paper is organized as follows. In Sect. 2 we briefly recall the basics
of shift systems and symbolic dynamics. The concepts and notation introduced
in this section (including the examples) will be used throughout. Order
patterns and forbidden root patterns, together with the outgrowth forbidden
patterns, are presented in Sect. 3. The structure of the outgrowth forbidden
patterns and their asymptotic growth with the length are discussed in Sect.
4. Finally, Sect. 5 and 6 are devoted to the structure of allowed patterns
and the existence of root forbidden patterns in one-sided and two-sided
shift systems, respectively. In the examples we present some interesting
by-products of the theoretical results.

\section{Shift systems and symbolic dynamics}

Let us start by recalling some basics of shift systems and symbolic
dynamics. We set $\mathbb{N}_{0}=\{0\}\cup \mathbb{N}=\{0,1,2,...\}$.

Fix $N\geq 2$ and consider the measurable space $(\Omega ,\mathcal{P}(\Omega
))$, where $\Omega =\{0,1,...,N-1\}$ and $\mathcal{P}(\Omega )$ is the
family of all subsets of $\Omega $. Let $(\Omega ^{\mathbb{N}_{0}},\mathcal{B%
})$ denote the product space $\Pi _{0}^{\infty }(\Omega ,\mathcal{P}(\Omega
))$, i.e., $\Omega ^{\mathbb{N}_{0}}$ is the space of (\textit{one-sided})
\textit{sequences} taking values on the `alphabet' $\Omega $,
\begin{equation*}
\Omega ^{\mathbb{N}_{0}}=\{\omega =(\omega _{n})_{n\in \mathbb{N}%
_{0}}:\omega _{n}\in \Omega \},
\end{equation*}%
and $\mathcal{B}$ is the sigma-algebra generated by the \textit{cylinder sets%
}
\begin{equation*}
C_{a_{0},...,a_{n}}=\{\omega \in \Omega ^{\mathbb{N}_{0}}:\omega
_{k}=a_{k},0\leq k\leq n\}.
\end{equation*}%
The topology generated by the cylinder sets makes $\Omega ^{\mathbb{N}_{0}}$
compact, perfect (i.e., it is closed and all its points are accumulation
points) and totally disconnected. Such topological spaces are sometimes
called Cantor sets. The elements of $\Omega $ are called \textit{symbols} or
\textit{letters}. Segments of symbols of length $L$, like $\omega _{k}\omega
_{k+1}\ldots \omega _{k+L-1}$, will be sometimes shortened $\omega
_{k}^{k+L-1}$.

Furthermore, let $\Sigma :\Omega ^{\mathbb{N}_{0}}\rightarrow \Omega ^{%
\mathbb{N}_{0}}$ denote the (one-sided) \textit{shift transformation}
defined as%
\begin{equation}
\Sigma :(\omega _{0},\omega _{1},\omega _{2},...)\mapsto (\omega _{1},\omega
_{2},\omega _{3},...).  \label{shift}
\end{equation}%
All probability measures on $(\Omega ^{\mathbb{N}_{0}},\mathcal{B})$ which
make $\Sigma $ a measure-preserving transformation are obtained in the
following way \cite{Walters}. For any $n\geq 0$ and $a_{i}\in \Omega $, $%
0\leq i\leq n$, let a real number $p_{n}(a_{0},...,a_{n})$ be given such
that (i) $p_{n}(a_{0},...,a_{n})\geq 0$, (ii) $\sum_{a_{0}\in \Omega
}p_{0}(a_{0})=1$, and (iii) $p_{n}(a_{0},...,a_{n})=\sum_{a_{n+1}\in \Omega
}p_{n+1}(a_{0},...,a_{n},a_{n+1})$. If we define now%
\begin{equation*}
m(C_{a_{0},...,a_{n}})=p_{n}(a_{0},...,a_{n}),
\end{equation*}%
then $m$ can be extended to a probability measure on $(\Omega ^{\mathbb{N}%
_{0}},\mathcal{B})$. The resulting dynamical system, $(\Omega ^{\mathbb{N}%
_{0}},\mathcal{B},m,\Sigma )$ is called the \textit{one-sided shift space}.

\medskip

\noindent \textbf{Example 1.} (a) Let $\mathbf{p}=(p_{0},p_{1},...,p_{N-1})$%
, $N\geq 2$, be a probability vector with non-zero entries (i.e., $p_{i}>0$
and $\sum\nolimits_{i=0}^{N-1}p_{i}=1$). Set $%
p_{n}(a_{0},a_{1},...,a_{n})=p_{a_{0}}p_{a_{1}}...p_{a_{n}}$. The resulting
measure-preserving shift transformation is called the one-sided $\mathbf{p}$-%
\textit{Bernoulli shift}.

(b) Let $\mathbf{p}=(p_{0},p_{1},...,p_{N-1})$ be a probability vector as in
(a) and $P=(p_{ij})_{0\leq i,j\leq N-1}$ an $N\times N$ stochastic matrix
(i.e., $p_{ij}\geq 0$ and $\sum\nolimits_{i,j=0}^{N-1}p_{ij}=1$) such that $%
\sum\nolimits_{i=0}^{N-1}p_{i}p_{ij}=p_{j}$. Set then $%
p_{n}(a_{0},a_{1},...,a_{n})=p_{a_{0}}p_{a_{0}a_{1}}p_{a_{1}a_{2}}...p_{a_{n-1}a_{n}}
$. The resulting measure-preserving shift transformation is called the
one-sided $(\mathbf{p},P)$-\textit{Markov shift}.

(c) Let $\mathbf{S}=(S_{n})_{n=0}^{\infty }$ be a discrete-time stochastic
process on a probability space $(X,\mathcal{F},\mu )$ started at time $n=0$
with finitely many outcomes $\{0,1,...,N-1\}=\Omega $. The realizations (or
\textquotedblleft sample paths\textquotedblright ) $\mathbf{S}%
(x)=(S_{0}(x),...,S_{n}(x),...)$ are viewed as elements of $\Omega ^{\mathbb{%
N}_{0}}$ endowed with the induced measure $p_{n}(a_{0},...,a_{n})=\mu
(\{x\in X:S_{0}(x)=a_{0},\ldots ,S_{n}(x)=a_{n}\})\equiv \Pr
\{S_{0}=a_{0},...,S_{n}=a_{n}\}$, the probability of the event $%
S_{0}=a_{0},...,S_{n}=a_{n}$. The resulting measure on $\Omega ^{\mathbb{N}%
_{0}}$ is shift invariant if the stochastic process $\mathbf{S}$ is
stationary. $\square $

\medskip

There are several metrics compatible with the topology of $\Omega ^{\mathbb{N%
}_{0}}$, the most popular being%
\begin{equation}
d_{K}(\omega ,\omega ^{\prime })=\sum\limits_{n=0}^{\infty }\frac{\delta
(\omega _{n},\omega _{n}^{\prime })}{K^{n}},  \label{metrick}
\end{equation}%
where $\delta (\omega _{n},\omega _{n}^{\prime })=1$ if $\omega _{n}\neq
\omega _{n}^{\prime }$, $\delta (\omega _{n},\omega _{n})=0$ and $K>2$.
Observe that given $\omega \in C_{a_{0},...,a_{n}}$, then $d_{K}(\omega
,\omega ^{\prime })<\frac{1}{K^{n}}$ if $\omega ^{\prime }\in
C_{a_{0},...,a_{n}}$ and $d_{K}(\omega ,\omega ^{\prime })\geq \frac{1}{K^{n}%
}$ if $\omega ^{\prime }\notin C_{a_{0},...,a_{n}}$, so that $%
C_{a_{0},...,a_{n}}=B_{d_{K}}(\omega ;\frac{1}{K^{n}})$, the open ball of
radius $K^{-n}$ and center $\omega $ in the metric space $(\Omega ^{\mathbb{N%
}_{0}},d_{K})$. Since the base of the measurable sets are open balls, we
conclude that $\mathcal{B}$ is the Borel sigma-algebra in the topology
defined by the metric (\ref{metrick}). Observe furthermore that every point
in $B_{d_{K}}(\omega ;\frac{1}{K^{n}})$ is a center, a property known from
non-Archimedean normed spaces (e.g., the rational numbers with p-adic norms).

Continuity will play a role below. Since $\Sigma
^{-1}C_{a_{0},...,a_{n}}=\cup _{a\in \Omega }C_{a,a_{0},...,a_{n}}$, $\Sigma
$ is continuous in $(\Omega ^{\mathbb{N}_{0}},d_{K})$, each point $\omega
\in \Omega ^{\mathbb{N}_{0}}$ having exactly $N$ preimages under $\Sigma $.
Regarding the forward dynamic, $\Sigma $ has $N$ fixed points: $\omega =(%
\bar{n})$, $0\leq n\leq N-1$, where \textit{the overbar denotes indefinite
repetition throughout}.

The corresponding (invertible) two-sided shift transformation on the
two-sided sequence (or \textit{bisequence}) space%
\begin{equation*}
\Omega ^{\mathbb{Z}}=\{\omega =(\omega _{n})_{n\in \mathbb{Z}}:\omega
_{n}\in \Omega \},
\end{equation*}%
is defined as $\Sigma :\omega \mapsto \omega ^{\prime }$ with $\omega
_{n}^{\prime }=\omega _{n+1}$, $n\in \mathbb{Z}$. The cylinder sets are
given now as $\{\omega \in \Omega ^{\mathbb{Z}}:\omega _{k}=a_{k},\left\vert
k\right\vert \leq n\}$ and%
\begin{equation*}
d_{K}(\omega ,\omega ^{\prime })=\sum\limits_{n\in \mathbb{Z}}\frac{\delta
(\omega _{n},\omega _{n}^{\prime })}{K^{\left\vert n\right\vert }},
\end{equation*}%
with $K>3$.

Let $T$ be a measure preserving map on a probability space $(X,\mathcal{F}%
,\mu )$ and $\alpha =\{A_{0},\ldots ,A_{N-1}\}$ be a generating partition of
the sigma-algebra $\mathcal{F}$ with respect to $T$, i.e., the subsets of
the form $A_{a_{0}}\cap T^{-1}A_{a_{1}}\cap ...\cap T^{-n}A_{a_{n}}$
generate $\mathcal{F}$. Assume moreover that for every sequence $%
(A_{a_{n}})_{n\in \mathbb{N}_{0}}$, the set $\cap _{n=0}^{\infty
}T^{-n}A_{a_{n}}$ contains at most one point of $X$; this assumption is
fulfilled by any positively expansive continuous map or expansive
homeomorphism on compact metric spaces (in particular, by the one-sided and
two-sided transformations we considered above) and implies that the coding
map $\Phi $ to be defined in (\ref{coding})-(\ref{coding2}) is one-to-one.
Define now on the cylinder sets of $\Omega ^{\mathbb{N}_{0}}$ the measure

\begin{equation*}
m_{T}(C_{a_{0},...,a_{n}})=\mu (A_{a_{0}}\cap T^{-1}A_{a_{1}}\cap ...\cap
T^{-n}A_{a_{n}}).
\end{equation*}

For $\omega \in \Omega ^{\mathbb{N}_{0}}$ define the \textit{coding map} $%
\Phi :X\rightarrow \Omega ^{\mathbb{N}_{0}}$ by%
\begin{equation}
\Phi (x)=(\omega _{0},...,\omega _{n},...),  \label{coding}
\end{equation}%
where%
\begin{equation}
\omega _{n}=a_{n}\in \Omega \;\;\text{\mbox{if}\ \ }T^{n}(x)\in
A_{a_{n}},n\geq 0.  \label{coding2}
\end{equation}%
Then $\Phi :(X,\mathcal{F},\mu )\rightarrow (\Omega ^{\mathbb{N}_{0}},%
\mathcal{B},m_{T})$ is measure-preserving (since, by definition, $\Phi
^{-1}(C_{a_{0},...,a_{n}})=A_{a_{0}}\cap T^{-1}A_{a_{1}}\cap ...\cap
T^{-n}A_{a_{n}}$) and, moreover,%
\begin{equation}
\Phi \circ T=\Sigma \circ \Phi ,  \label{iso}
\end{equation}%
i.e., $T$ and $\Sigma $ are isomorphic and, hence, $(X,\mathcal{F},\mu ,T)$
and $(\Omega ^{\mathbb{N}_{0}},\mathcal{B},m_{T},\Sigma )$ are dynamically
equivalent.

One interesting consequence of this construction is that the coded orbits of
$T$ contain any arbitrary pattern. Indeed, given any $N$-symbol \textit{%
pattern} of length $L\geq 1$, $a_{0}^{L-1}:=a_{0}a_{1}...a_{L-1}$ with
symbols $a_{n}\in \{0,1,...,N-1\}$, choose
\begin{equation*}
x_{0}\in \bigcap\limits_{n=0}^{L-1}T^{-n}A_{a_{n}}.
\end{equation*}%
Then $\Phi (x_{0})\in C_{a_{0},...,a_{L-1}}$ and this for any $L\geq 1$.
Letting $L\rightarrow \infty $, we conclude that the coding map $\Phi $
associates to each orbit $orb(x)=\{T^{n}(x):n\geq 0\}$ a unique, infinitely
long pattern of symbols from $\{0,1,...,N-1\}$, namely, $\Phi (x)$, for
almost all $x\in X$.

In the special case of invertible maps $T:X\rightarrow X$, both $T$ and $%
T^{-1}$ are measurable and all the above generalizes to two-sided sequences.

\medskip

\noindent \textbf{Example 2.} As a standard example (that it is going to be
our workhorse), take $X=[0,1],$ $\mathcal{F}$ the Borel sigma-algebra
restricted to $[0,1]$, $d\mu =\frac{1}{\pi \sqrt{x(1-x)}}dx$, $f(x)=4x(1-x)$%
, the \textit{logistic map}, and $\alpha =\{A_{0}=[0,\frac{1}{2}),A_{1}=[%
\frac{1}{2},1]\}$ (it is irrelevant whether the midpoint $\frac{1}{2}$
belongs to the left or to the right partition element). Then $\Phi (\frac{1}{%
4})=(0,\bar{1})$, $\Phi (\frac{1}{2})=(1,1,\bar{0})$ and $\Phi (\frac{3}{4}%
)=(\bar{1})$. Observe for further reference that $\Phi (\frac{1}{4})<\Phi (%
\frac{1}{2})<\Phi (\frac{3}{4})$, where $<$ stands for the lexicographical
order of $\{0,1\}^{\mathbb{N}_{0}}$, but, e.g., $\Phi (\frac{1}{2})>\Phi
(1)=(1,\bar{0})$, hence the coding map $\Phi :[0,1]\rightarrow \{0,1\}^{%
\mathbb{N}_{0}}$ does not preserve the order structure. The fixed points of $%
f$ are $0=\Phi ^{-1}((\bar{0}))$ and $\frac{3}{4}=\Phi ^{-1}((\bar{1}))$. $%
\square$

\medskip

\section{\noindent \textbf{Forbidden order patterns}}

In the previous Section, we saw that the symbolic dynamics of maps defines
any symbol pattern of any length, under rather general assumptions. In this
Section we will see that the situation is not quite the same when
considering order patterns.

Let $(X,<)$ be a totally ordered set and $T:X\rightarrow X$ a map. Given $%
x\in X$, the orbit of $x$ is the set $\{T^{n}(x):n\in \mathbb{N}_{0}\}$,
where $T^{0}(x)\equiv x$ and $T^{n}(x)\equiv T(T^{n-1}(x))$. If $x$ is not a
periodic point of period less than $L\geq 2$, we can then associate with $x$
an order pattern of length $L$, as follows. We say that $x$ \textit{defines
the order pattern} $\pi =\pi (x)=[\pi _{0},...,\pi _{L-1}]$, where $\{\pi
_{0},...,\pi _{L-1}\}$ is a permutation of $\{0,1,\ldots ,L-1\}$, if
\begin{equation*}
T^{\pi _{0}}(x)<T^{\pi _{1}}(x)<...<T^{\pi _{L-1}}(x).
\end{equation*}%
Alternatively, we say that $x$ \textit{is of type} $\pi $ or that $\pi $
\textit{is realized by} $x$. Thus, $\pi $ is just a permutation on $%
\{0,1,...,L-1\}$, given by $0\mapsto \pi _{0},...,L-1\mapsto \pi _{L-1}$,
that encapsulates the order of the points $x_{n}=T^{n}(x)$, $0\leq n\leq L-1$%
. The set of order patterns of length $L$ or, equivalently, the set of
permutations on $\{0,1,...,L-1\}$ will be denoted by $\mathcal{S}_{L}$.
Furthermore set%
\begin{equation*}
P_{\pi }=\{x\in X:x\text{ defines }\pi \in \mathcal{S}_{L}\}.
\end{equation*}

A plain difference between symbol patterns and order patterns of length $L$
is their cardinality: the former grow exponentially with $L$ (exactly as $%
N^{L}$, where $N$ is the number of symbols$)$ while the latter do
super-exponentially. Specifically,%
\begin{equation}
\left\vert \mathcal{S}_{L}\right\vert =L!\propto e^{L(\ln L-1)+(1/2)\ln 2\pi
L}  \label{Stirling}
\end{equation}%
(Stirling's formula), where, as usual, $\left\vert \cdot \right\vert $
denotes cardinality and $\propto $ means \textquotedblleft
asymptotically\textquotedblright . Although one can construct functions
whose orbits realize any possible order pattern (see below), numerical
simulations support the conjecture that order patterns, like symbol
patterns, grow only exponentially for `well-behaved' functions \cite{Bandt}.
In fact, if $I$ is a closed interval of $\mathbb{R}$ and $f:I\rightarrow I$
is \textit{piecewise monotone} (i.e., there is a finite partition of $I$
into intervals such that $f$ is continuous and strictly monotone on each of
those intervals), then one can prove \cite{Bandt} that
\begin{equation}
\left\vert \{\pi \in \mathcal{S}_{L}:P_{\pi }\neq \varnothing \}\right\vert
\propto e^{Lh_{top}(f)},  \label{TopEntropy2}
\end{equation}%
where $h_{top}(f)$ is the topological entropy of $f$. From (\ref{Stirling})
and (\ref{TopEntropy2}) we conclude:

\medskip

\noindent \textbf{Proposition 1. }If\textit{\ }$f$ is a piecewise monotone
map on a closed interval\textit{\ }$I\subset \mathbb{R}$, then there are%
\textit{\ }$\pi \in \mathcal{S}_{L}$, $L\geq 2$,\ such that\textit{\ }$%
P_{\pi }=\varnothing $.

\medskip

Order patterns that do not appear in any orbit of $f$ are called \textit{%
forbidden patterns}, at variance with the \textit{allowed patterns}, for
which there are intervals of points that realize them.%

\begin{figure}[hbt]
\begin{center}
\DeclareGraphicsExtensions{eps} \includegraphics[width=1.0%
\textwidth]{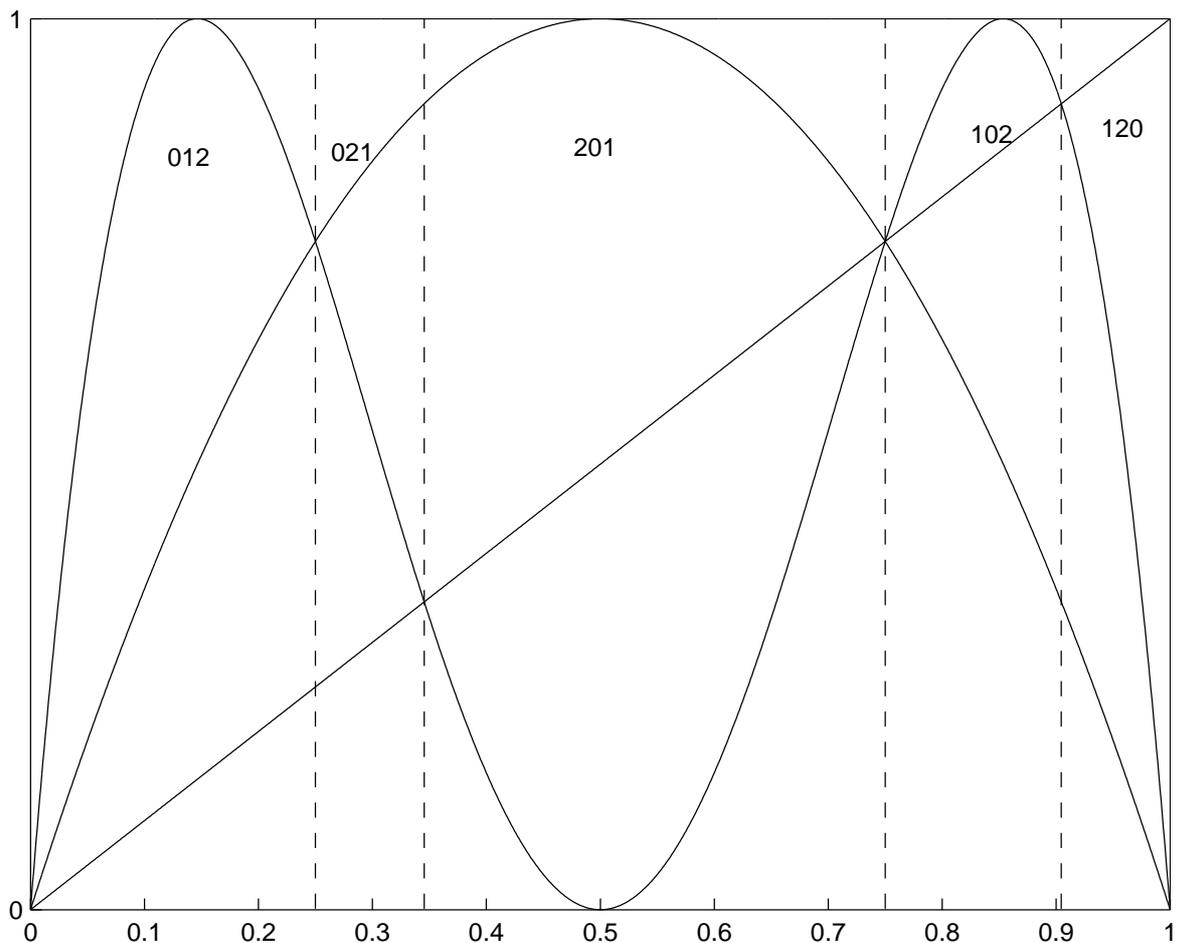}
\end{center}
\par
\caption{The sets $P_{\pi }$, $\pi \in \sigma _{3}$, are graphically
obtained
by raising vertical lines at the crossing points of the curves $y=x$, $%
y=f(x) $, and $y=f^{2}(x)$. The three digits on the top are
shorthand for
order patterns (e.g., $012$ stands for $[0,1,2]$). We see that $%
P_{[2,1,0]}=\varnothing $.}
 \label{ndecrease}
\end{figure}

\noindent \textbf{Example 3.} As a simple illustration borrowed from \cite%
{Amigo2}, consider again the logistic map. For $L=2$ we have%
\begin{equation*}
P_{[0,1]}=\left( 0,\tfrac{3}{4}\right) ,\;\;P_{[1,0]}=\left( \tfrac{3}{4}%
,1\right) .
\end{equation*}%
Observe that the endpoints of $P_{\pi }$ are period-1 (i.e., fixed) points ($%
0$ and $\tfrac{3}{4}$) or preimages of them ($f(1)=0$). But already for $L=3$
($f^{2}(x)=-64x^{4}+128x^{3}-80x^{2}+16x$) there are permutations that are
not realized (see Figure~\ref{ndecrease}):%
\begin{equation*}
\begin{array}{lll}
P_{[0,1,2]}=\left( 0,\frac{1}{4}\right) , & P_{[0,2,1]}=\left( \frac{1}{4},%
\frac{5-\sqrt{5}}{8}\right) , & P_{[2,0,1]}=\left( \frac{5-\sqrt{5}}{8},%
\frac{3}{4}\right) , \\
P_{[1,0,2]}=\left( \frac{3}{4},\frac{5+\sqrt{5}}{8}\right) , &
P_{[1,2,0]}=\left( \frac{5+\sqrt{5}}{8},1\right) , & P_{[2,1,0]}=\varnothing
.%
\end{array}%
\end{equation*}%
In going from $\pi \in \mathcal{S}_{2}$ to $\pi \in \mathcal{S}_{3}$, we see
that $P_{[0,1]}$ splits into the subintervals $P_{[0,1,2]}$, $P_{[0,2,1]}$
and $P_{[2,0,1]}$ at the eventually periodic point $\frac{1}{4}$ (preimage
of $\frac{3}{4}$) and at the period-2 point $\frac{5-\sqrt{5}}{8}$.
Likewise, $P_{[1,0]}$ splits into $P_{[1,0,2]}$ and $P_{[1,2,0]}$ at the
period-2 point $\frac{5+\sqrt{5}}{8}$.

From a different perspective, as we move rightward in
Figure~\ref{ndecrease} from the
neighborhood of $0$, where $x<f(x)<f^{2}(x)$, the curves $y=f(x)$ and $%
y=f^{2}(x)$ cross at $x=\frac{1}{4}$, what causes the first swap: $[0,1,2]$
transforms to $[0,2,1]$. In general, the crossings at $x=\frac{1}{4},\frac{5-%
\sqrt{5}}{8}$ and $\frac{5+\sqrt{5}}{8}$ between $f^{\pi (i)}$ and $f^{\pi
(i+1)}$ causes the exchange of $\pi (i)$ and $\pi (i+1)$ in the pre-crossing
pattern. At $x=\frac{3}{4}$ all three curves cross and $[2,0,1]$ goes over
to $[1,0,2]$.

The absence of $\pi =[2,1,0]$ triggers, in turn, an avalanche of longer
missing patterns. To begin with, the pattern $[\ast ,2,\ast ,1,\ast ,0,\ast
] $ (where the wildcard $\ast $ stands eventually for any other entries of
the pattern) cannot be realized by any $x\in \lbrack 0,1]$ since the
inequality
\begin{equation}
f^{2}(x)<f(x)<x  \label{inequality}
\end{equation}%
cannot occur. By the same token, the patterns $[\ast ,3,\ast ,2,\ast ,1,\ast
]$, $[\ast ,4,\ast ,3,\ast ,2,\ast ]$, and, more generally, $[\ast ,n+2,\ast
,n+1,\ast ,n,\ast ]\in \mathcal{S}_{L}$, $0\leq n\leq L-3$, cannot be
realized either for the same reason (substitute $x$ by $f^{n}(x)$ in (\ref%
{inequality})). $\square$

\medskip

The same follows for the \textit{tent map} $\Lambda :[0,1]\rightarrow
\lbrack 0,1]$,
\begin{equation}
\Lambda (x)=\left\{
\begin{array}{cl}
2x & \;\;0\leq x\leq \frac{1}{2} \\
2-2x & \;\;\frac{1}{2}\leq x\leq 1%
\end{array}%
\right. .  \label{tent}
\end{equation}%
In fact, if $\lambda $ is the Lebesgue measure, $d\mu =\frac{1}{\pi \sqrt{%
x(1-x)}}dx$ is (as in Example 2) the invariant measure of the logistic map $%
f(x)=4x(1-x)$, and $\phi :([0,1],\lambda )\rightarrow ([0,1],\mu )$ is the
measure preserving isomorphism given by%
\begin{equation}
\phi (x)=\sin ^{2}(\tfrac{\pi }{2}x),  \label{phi}
\end{equation}%
then the dynamical systems $([0,1],\mathcal{B},\lambda ,\Lambda )$ and $%
([0,1],\mathcal{B},\mu ,f)$, where $\mathcal{B}$ is the Borel sigma-algebra
restricted to the interval $[0,1]$, are \textit{isomorphic} (or conjugate)
by means of $\phi $, i.e., $f\circ \phi =\phi \circ \Lambda $. Since,
moreover, $\phi $ is strictly increasing, forbidden patterns for $f$
correspond to forbidden patterns for $\Lambda $ in a one-to-one way.

From the last paragraph it should be clear that isomorphic dynamical systems
need not have the same forbidden patterns: the isomorphism ($\phi $ above)
must also preserve the linear order of both spaces (supposing both spaces
are linearly ordered), and this will be in general not the case. For
example, the $\lambda $-preserving \textit{shift map} $S_{2}:x\mapsto 2x(%
\mbox{mod 1})$, $0\leq x\leq 1$, has no forbidden patterns of length $3$,
although it is isomorphic to the logistic and tent maps (the isomorphism
with $f$ is proved via the semi-conjugacy $\varphi :([0,1],\lambda
)\rightarrow ([0,1],\mu )$, $\varphi (x)=\sin ^{2}\pi x$, which does not
preserve order on account of being increasing on $(0,\frac{1}{2})$ and
decreasing on $(\frac{1}{2},1)$). The same happens with the logistic map and
the $(\frac{1}{2},\frac{1}{2})$-Bernoulli shift, a model for tossing of a
fair coin, because, as we saw in Example 2, the corresponding isomorphy
(actually, the coding map) $\Phi :[0,1]\rightarrow \{0,1\}^{\mathbb{N}_{0}}$
is not order-preserving.

Two isomorphic dynamical systems, whose phase spaces are linearly ordered,
are called \textit{order-isomorphic} if the isomorphism between them is also
an order-isomorphism (i.e., it also preserves the order structure). It is
obvious that two order-isomorphic systems (like those defined by the
logistic and the tent map) have the same order patterns.

\medskip

\noindent \textbf{Proposition 2.} Given $X_{1},X_{2}\subset \mathbb{R}$
endowed with the standard Borel sigma-algebra $\mathcal{B}$, suppose that
the dynamical systems $(X_{1},\mathcal{B},\mu _{1},f_{1})$ and $(X_{2},%
\mathcal{B},\mu _{2},f_{2})$ are isomorphic via a continuous map $\phi
:X_{1}\rightarrow X_{2}$. If $f_{1}$ is topologically transitive and, for
all $x\in X_{1}$, both $x$ and $\phi (x)$ define the same order patterns,
then $\phi $ is order-preserving.

\medskip

\noindent \textit{Proof}. We want to prove that if $x,x^{\prime }\in X_{1}$
and $x<x^{\prime }$, then $\phi (x)<\phi (x^{\prime }).$ Because of
continuity, for all $\varepsilon >0$ there exists $0<\delta <\frac{x^{\prime
}-x}{2}$ such that $\left\vert y-x\right\vert <\delta \Rightarrow \left\vert
\phi (y)-\phi (x)\right\vert <\frac{\varepsilon }{2}$ and $\left\vert
y^{\prime }-x^{\prime }\right\vert <\delta $ $\Rightarrow \left\vert \phi
(y^{\prime })-\phi (x^{\prime })\right\vert <\frac{\varepsilon }{2}$. On the
other hand, transitiveness implies that, given $x$, $x^{\prime }$ and $%
\delta $ as above, there exists $x_{0}\in X_{1}$, $N=N(x,\delta )$ and $%
N^{\prime }=N^{\prime }(x^{\prime },\delta )$ such that $\left\vert
f_{1}^{N}(x_{0})-x\right\vert <\delta $ and $\left\vert f_{1}^{N^{\prime
}}(x_{0})-x^{\prime }\right\vert <\delta $. Thus $%
f_{1}^{N}(x_{0})<f_{1}^{N^{\prime }}(x_{0})$ and, by assumption, $\phi \circ
f_{1}^{N}(x_{0})=f_{1}^{N}(\phi (x_{0}))<f_{2}^{N^{\prime }}(\phi
(x_{0}))=\phi \circ f_{1}^{N^{\prime }}(x_{0})$ holds. Choose now $%
y=f_{1}^{N}(x_{0})$, $y^{\prime }=f_{1}^{N^{\prime }}(x_{0})$ to deduce%
\begin{equation*}
\phi (y)<\phi (y^{\prime })\leq \phi (x^{\prime })+\left\vert \phi
(y^{\prime })-\phi (x^{\prime })\right\vert \leq \phi (x^{\prime })+\frac{%
\varepsilon }{2},
\end{equation*}%
where $\varepsilon $ is arbitrary. If we choose now $\varepsilon <\frac{%
\left\vert \phi (x)-\phi (x^{\prime })\right\vert }{2}$, then it follows $%
\phi (x)<\phi (x^{\prime })$, since $\left\vert \phi (y)-\phi (x)\right\vert
<\frac{\varepsilon }{2}$. $\square $

\medskip

Finally, observe that the setting we are considering is more general than
the setting of Kneading Theory since our functions need not be continuous
(but only piecewise-continuous). Under some assumptions \cite{Melo}, the
kneading invariants completely characterize the order-isomorphy of
continuous maps.

\medskip

\section{\noindent \textbf{Outgrowth forbidden patterns }}

According to Proposition 1, for every piecewise monotone interval map on $%
\mathbb{R}$, $f:I\rightarrow I$, there exist $\pi \in \mathcal{S}_{L}$, $%
L\geq 2$, which cannot occur in any orbit. We call them \textit{forbidden
patterns} for $f$ and recall how their absence pervades all longer patterns
in form of \textit{outgrowth forbidden patterns }(see Example 3). Since $\pi
=[\pi _{0},...,\pi _{L-1}]$ is forbidden for $f$, then the $2(L+1)$ patterns
of length $L+1$,
\begin{eqnarray*}
&&[L,\pi _{0},...,\pi _{L-1}],[\pi _{0},L,\pi _{1},...,\pi _{L-1}],...,[\pi
_{0},...,\pi _{L-1},L], \\
&&[0,\pi _{0}+1,...,\pi _{L-1}+1],[\pi _{0}+1,0,\pi _{1}+1,...,\pi
_{L-1}+1],...,[\pi _{0}+1,...,\pi _{L-1}+1,0],
\end{eqnarray*}%
are also forbidden for $f$. Assume for the time being that all these
forbidden patterns belonging to the \textquotedblleft first
generation\textquotedblright\ are all different. Then, proceeding similarly
as before, we would find%
\begin{equation*}
2(L+1)\times 2(L+2)=2^{2}(L+1)(L+2)
\end{equation*}%
forbidden patterns of length $L+2$ in the second generation and, in general,%
\begin{equation*}
2^{m}(L+1)...(L+m)=2^{n}\frac{(L+m)!}{L!}
\end{equation*}%
forbidden patterns of length $L+m$ in the $m$th generation, provided that
all forbidden patterns up to (and including) the $m$th generation are
different. Observe that all these forbidden patterns generated by $\pi $
have the form
\begin{equation}
\lbrack \ast ,\pi _{0}+n,\ast ,\pi _{1}+n,\ast ,...,\ast ,\pi _{L-1}+n,\ast
]\in \mathcal{S}_{N}  \label{piele}
\end{equation}%
with $n=0,1,...,N-L$, where $N-L\geq 1$ is the number of wildcards $\ast \in
\{0,1,...,n-1,L+n,...,N-1\}$ (with $\ast \in \{L,...,N-1\}$ if $n=0$ and $%
\ast \in \{0,...,N-L-1\}$ if $n=N-L$).

A better upper bound on the number of outgrowth forbidden patterns of length
$N$ of $\pi $ is obtained using the following reasoning. For fixed $n$, the
number of outgrowth patterns of $\pi $ of the form (\ref{piele}) is $%
N!/(N-L)!$. This is because out of all possible permutations of the numbers $%
\{0,1,\ldots ,N-1\}$, we only count those that have the entries $\{\pi
_{0}+n,\pi _{1}+n,\ldots ,\pi _{L-1}+n\}$ in the required order. Next, note
that we have $N-L+1$ choices for the value of $n$. Each choice generates a
set of $N!/(N-L)!$ outgrowth patterns. These sets are not necessarily
disjoint, but an upper bound on the size of their union, i.e., the set of
all outgrowth forbidden patterns of length $N$ of $\pi $, is given by
\begin{equation*}
(N-L+1)\frac{N!}{(N-L)!}.
\end{equation*}

A weak form of the converse holds also true: if $[L,\pi _{0},...,\pi _{L-1}]$%
, $[\pi _{0},L,...,\pi _{L-1}]$, $...$, $[\pi _{0},...,\pi _{L_{0}-1},L]\in
\mathcal{S}_{L+1}$ are forbidden, then $[\pi _{0},...,\pi _{L-1}]\in
\mathcal{S}_{L}$ is also forbidden.

Forbidden patterns that are not outgrowth patterns of other forbidden
patterns of shorter length are called \textit{forbidden root patterns} since
they can be viewed as the root of the tree of forbidden patterns spanned by
the outgrowth patterns they generate, branching taking place when going from
one length (or generation) to the next.%

\begin{figure}[htb]
\begin{center}
\DeclareGraphicsExtensions{eps} \includegraphics[width=1.0%
\textwidth]{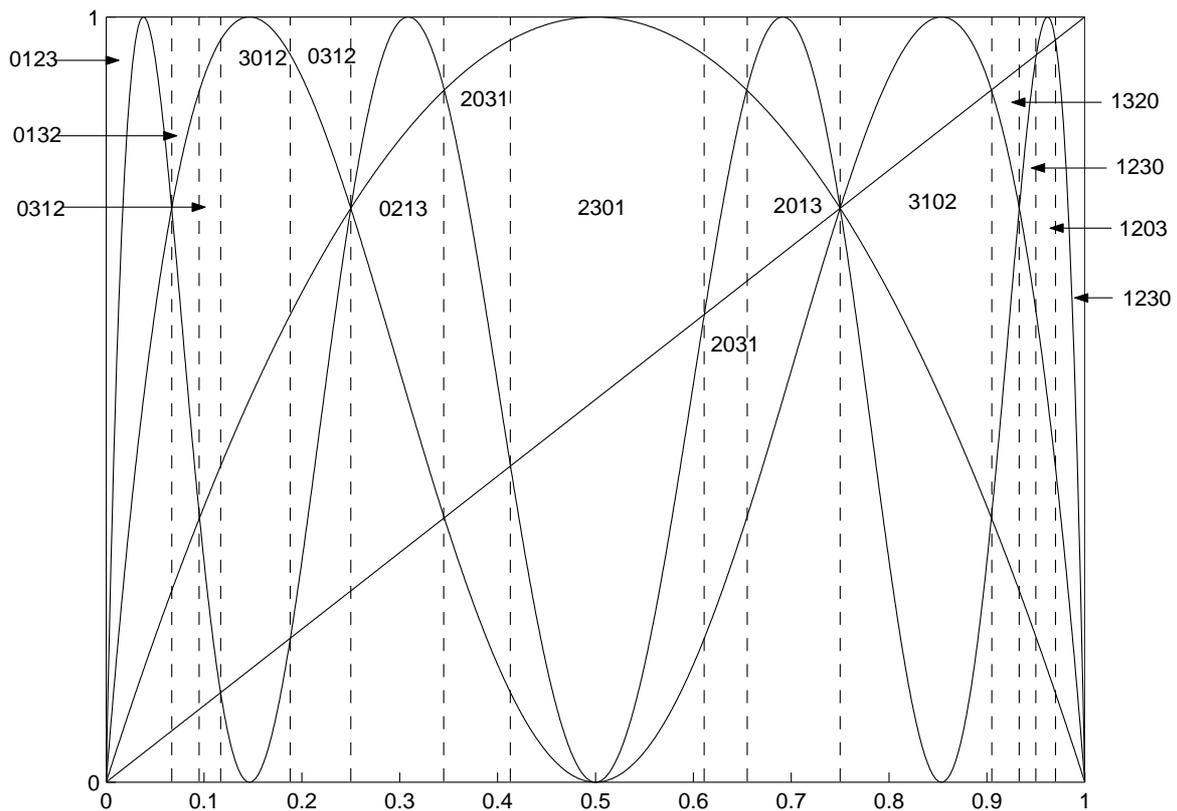}
\end{center}
\par
\caption{The twelve allowed order patterns of length 4 for the
logistic map. Note the two components of $P_{[0,3,1,2]}$,
$P_{[2,0,3,1]}$ and $P_{[1,2,3,0]}$.} \label{Fig2}
\end{figure}

\noindent \textbf{Example 4.} If $f$ is the logistic map, then%
\begin{equation*}
f^{3}(x)=-16\,384x^{8}+65\,536x^{7}-106\,496x^{6}+90\,112\allowbreak
x^{5}-42\,240x^{4}+10\,752x^{3}-1344x^{2}+64x.
\end{equation*}%
In Figure~\ref{Fig2}, which is Figure~\ref{ndecrease} with the curve
$y=f^{3}(x)$ superimposed, we can see the 12 allowed patterns of
length 4 of the logistic map. Since there are 24 possible patterns
of length 4, we conclude that 12 of them are forbidden. The
outgrowth patterns of $[2,1,0]$, the only forbidden pattern of
length 3,
are (see (\ref{piele})):%
\begin{equation*}
\begin{array}{cl}
(n=0) & [3,2,1,0],[2,3,1,0],[2,1,3,0],[2,1,0,3] \\
(n=1) & [0,3,2,1],[3,0,2,1],[3,2,0,1],[3,2,1,0]%
\end{array}%
.
\end{equation*}%
Observe that the pattern $[3,2,1,0]$ is repeated. Therefore, the remaining
five forbidden patterns of length 4 are root patterns.

In Figure~\ref{Fig2} one can also follow the first two splittings of the intervals $%
P_{\pi }$:
\begin{eqnarray*}
P_{[0,1]} &\rightarrow &\left\{
\begin{array}{l}
P_{[0,1,2]}\rightarrow
P_{[0,1,2,3]},P_{[0,1,3,2]},P_{[0,3,1,2]},P_{[3,0,1,2]} \\
P_{[0,2,1]}\rightarrow P_{[0,2,1,3]} \\
P_{[2,0,1]}\rightarrow P_{[2,0,1,3]},P_{[2,0,3,1]},P_{[2,3,0,1]}%
\end{array}%
\right. , \\
P_{[1,0]} &\rightarrow &\left\{
\begin{array}{l}
P_{[1,0,2]}\rightarrow P_{[3,1,0,2]} \\
P_{[1,2,0]}\rightarrow P_{[1,2,0,3]},P_{[1,2,3,0]},P_{[1,3,2,0]}%
\end{array}%
\right. .
\end{eqnarray*}%
The splitting of the intervals $P_{\pi }$ can be understood in terms of
periodic points and their preimages. Thus, the splitting of $P_{[0,1]}$ is
due to the points $\frac{1}{4}$ (first preimage of the period-1 point $\frac{%
3}{4}$) and $\frac{5-\sqrt{5}}{8}$ (a period-2 point); the second period-2
point, $\frac{5-\sqrt{5}}{8}$, is responsible for the splitting of $%
P_{[1,0]} $. On the contrary, $P_{[0,2,1]}$ and $P_{[1,0,2]}$ do not split
because they contain neither period-3 point nor first preimages of period-2
points nor second preimages of fixed points. $\square $

\medskip

Given the permutation $\sigma \in \mathcal{S}_{N},$ we say that $\sigma $
contains the \textit{consecutive pattern} $\tau =[\tau _{0},\tau _{1},\ldots
,\tau _{L-1}]\in \mathcal{S}_{L}$, $L<N$, if it contains a consecutive
subsequence order-isomorphic to $\tau $. Alternatively, we say that $\sigma $
avoids the \textit{consecutive pattern} $\tau $ if it contains no
consecutive subsequence order-isomorphic to $\tau $ \cite{EliI}.

Suppose now $\sigma \in \mathcal{S}_{N}$, $\pi \in \mathcal{S}_{L}$, $L<N$,
and%
\begin{equation*}
\begin{array}{llll}
\pi (p_{0})=0, & \pi (p_{1})=1, & \ldots , & \pi (p_{L-1})=L-1, \\
\sigma (s_{0})=n & \sigma (s_{1})=1+n, & \ldots , & \sigma (s_{L-1})=L-1+n,%
\end{array}%
\end{equation*}%
with $n\in \{0,1,...,N-L\}$. Then, the sequences $p_{0},p_{1},\ldots
,p_{L-1} $ and $s_{0},s_{1},\ldots ,s_{L-1}$ are consecutive subsequences of
$\pi ^{-1}$ and $\sigma ^{-1}$ (starting at positions $0$ and $n$),
respectively. If, moreover, $\sigma $ is an outgrowth pattern of $\pi $ (see
(\ref{piele})), then $s_{0},s_{1},\ldots ,s_{L-1}$ is order-isomorphic to $%
p_{0},p_{1},\ldots ,p_{L-1}$. It follows that $\sigma $ $\in \mathcal{S}_{N}$
is an outgrowth pattern of $\pi =[\pi _{0},\ldots ,\pi _{L-1}]$ if $\sigma
^{-1}$ contains $\pi ^{-1}$ as a consecutive subsequence. Hence, the allowed
patterns for $f$ are the permutations that avoid all such consecutive
subsequences for every forbidden root pattern of $f$.

\medskip

\noindent \textbf{Example 5.} Take $\pi =[2,0,1]$ to be a forbidden pattern
for a certain function $f$. Then $\sigma =[4,2,1,5,3,0]$ is an outgrowth
pattern of $\pi $ because it contains the subsequence $4,2,3$ ($n=2$).
Equivalently, $\sigma ^{-1}=[5,2,1,4,0,3]$ contains the consecutive pattern $%
1,4,0$ (starting at location $\sigma _{2}^{-1})$, which is order-isomorphic
to $\pi ^{-1}=[1,2,0]$. $\square $

\medskip

Let out$(\pi )$ denote the family of outgrowth patterns of the forbidden
pattern $\pi $,%
\begin{equation*}
\text{out}_{N}(\pi )=\text{out}(\pi )\cap \mathcal{S}_{N}=\{\sigma \in
\mathcal{S}_{N}:\sigma ^{-1}\text{ contains }\pi ^{-1}\text{ as a
consecutive pattern}\},
\end{equation*}%
and%
\begin{equation*}
\text{avoid}_{N}(\pi )=\mathcal{S}_{N}\backslash \text{out}_{N}(\pi
)=\{\sigma \in \mathcal{S}_{N}:\sigma ^{-1}\text{ avoids }\pi ^{-1}\text{ as
a consecutive pattern}\}.
\end{equation*}%
where $\backslash $\ stands for set difference. The fact that some of the
outgrowth patterns of a given length will be the same and that this depends
on $\pi $, makes the analytical calculation of $\left\vert \text{out}%
_{N}(\pi )\right\vert $ extremely complicated. Yet, from \cite{EliI} we know
that there are constants $0<c,d<1$ such that%
\begin{equation*}
c^{N}N!<\left\vert \text{avoid}_{N}(\pi )\right\vert <d^{N}N!
\end{equation*}%
(for the first inequality, $L\geq 3$ is needed). This implies that%
\begin{equation}
(1-d^{N})N!<\left\vert \text{out}_{N}(\pi )\right\vert <(1-c^{N})N!.
\label{outN}
\end{equation}

This factorial growth with $N$ can be exploited in practical applications to
tell random from deterministic time series with, in principle, arbitrarily
high probability. As said in the Introduction, these practical aspects are
beyond the scope of this paper, but let us bring up here the following,
related point. In the case of real (hence, \textit{finite}) randomly
generated sequences, a given order pattern $\pi \in \mathcal{S}_{L}$ can be
missing with nonvanishing probability. We call \textit{false forbidden
patterns }such\textit{\ }missing order patterns in finite random sequences
without constraints, to distinguish them from the `true' forbidden patterns
of deterministic (finite or infinite) sequences. True and false forbidden
patterns of self maps on one-dimensional intervals have been studied in \cite%
{Amigo3}.

\section{\noindent \textbf{Order patterns and one-sided shifts}}

The general study of order patterns and forbidden patterns is quite
difficult. Analytical results seem to be only feasible for particular maps.
In this and next sections we will consider the one- and two-sided shifts
since, owing to their simple structure, they can be analyzed with greater
detail. As we saw in Sect. 2, shifts are continuous maps (automorphisms if
two-sided) on compact metric spaces $(\{0,1,$ $...,N-1\}^{\mathbb{N}%
_{0}},d_{K})$ (resp., $(\{0,1,$ $...,N-1\}^{\mathbb{Z}},d_{K})$) that can be
lexicographically ordered:
\begin{equation*}
\omega <\omega ^{\prime }\;\Leftrightarrow \left\{
\begin{array}{l}
\;\omega _{0}<\omega _{0}^{\prime } \\
\text{\mbox{or}} \\
\omega _{0}=\omega _{0}^{\prime },...,\text{ }\omega _{n-1}=\omega
_{n-1}^{\prime }\text{ \mbox{and} }\omega _{n}<\omega _{n}^{\prime }\;(n\geq
1)%
\end{array}%
\right. \text{ ,}
\end{equation*}%
If $\overline{\mathcal{N}}$ denotes the countable, dense and $\Sigma $%
-invariant set of $\omega $ eventually terminating in an infinite string of $%
(N-1)$s except the sequence $(\overline{N-1})$, then the map $\psi
:\{0,1,...,N-1\}^{\mathbb{N}_{0}}\backslash \overline{\mathcal{N}}%
\rightarrow \lbrack 0,1]$ defined by

\begin{equation}
\psi :(\omega _{n})_{n\in \mathbb{N}_{0}}\mapsto \sum\limits_{n=0}^{\infty
}\omega _{n}N^{-(n+1)}.  \label{iso2}
\end{equation}%
is one-to-one and \textit{order-preserving}; moreover, $\psi ^{-1}$ is also
order-preserving. As a matter of fact, the lexicographical order in $\{0,1,$
$...,N-1\}^{\mathbb{N}_{0}}\backslash \overline{\mathcal{N}}$ corresponds
via $\psi $ to the standard order (induced by the positive numbers) in the
interval $[0,1]$. Although not important for our purposes, let us point out
that $\psi $ is continuous, but $\psi ^{-1}$ is not. Since the map
\begin{equation}
S_{N}=\psi \circ \Sigma \circ \psi ^{-1}:[0,1]\rightarrow \lbrack 0,1],
\label{sawtooth}
\end{equation}%
where $\Sigma $ is the shift on $N$ symbols, is piecewise linear and $%
\overline{\mathcal{N}}$ is dense, it follows (Proposition 1) that $\Sigma $
will have forbidden order patterns (although $\Sigma $ has no forbidden
\textit{symbol} pattern, see Sect. 2). In particular, if $\Sigma $ is the $(%
\frac{1}{N},...,\frac{1}{N})$-Bernoulli shift, then $S_{N}$ is the
Lebesgue-measure preserving \textit{sawtooth map }$S_{N}:x\mapsto Nx(%
\mbox{mod 1})$. Observe that only sequences that are not eventually periodic
define order patterns of any length.

What is the structure of the allowed order patterns? It is easy to convince
oneself (see Example 6 below) that, given $\omega =(\omega _{0},...,\omega
_{L-1},...)\in \{0,1,...,N-1\}^{\mathbb{N}_{0}}$ of type $\pi \in \mathcal{S}%
_{L}$, $\pi $ can be decomposed into, in general, $N$ blocks,%
\begin{equation}
\lbrack \pi _{0},...,\pi _{k_{0}-1};\pi _{k_{0}},...,\pi
_{k_{0}+k_{1}-1};...;\pi _{k_{0}+...+k_{N-2}},...,\pi
_{k_{0}+...+k_{N-2}+k_{N-1}-1}],  \label{allowed}
\end{equation}%
the semicolons separating the different blocks, where $k_{n}\geq 0$, $0\leq
n\leq N-1$, is the number of symbols $n\in \{0,1,\ldots ,N-1\}$ in $\omega
_{0}^{L-1}$ ($k_{n}=0$ if none, with the corresponding block missing) and $%
k_{0}+...+k_{N-1}=L$. Moreover:

(R1) The first (leftmost) block, $\pi _{0},...,\pi _{k_{0}-1}$, contains the
locations of the $0$s in $\omega _{0}^{L-1}$. Each $0$-run (i.e., a segment
of two or more consecutive $0$s contained in or intersected by $\omega
_{0}^{L-1}$), if any, contributes an \textit{increasing} subsequence $\pi
_{i},\pi _{i}+1,\pi _{i}+2,\ldots $ (as long as the $0$-run), which is
possibly intertwined with other entries of this block.

(R2) The last (rightmost) block, $\pi _{k_{0}+...+k_{N-2}},...,\pi
_{k_{0}+...+k_{N-2}+k_{N-1}-1}$, contains the locations of the $(N-1)$s in $%
\omega _{0}^{L-1}$. Each $(N-1)$-run contained in or intersected by $\omega
_{0}^{L-1}$, if any, contributes a \textit{decreasing} subsequence $\pi
_{k_{0}+...+k_{N-2}+i},\pi _{k_{0}+...+k_{N-2}+i}-1,\ldots $ (as long as the
$(N-1)$-run), which is possibly intertwined with other entries of this block.

(R3) Every intermediate block, $\pi _{k_{0}+...+k_{j-1}},...,\pi
_{k_{0}+...+k_{j-1}+k_{j}-1}$, $1\leq j\leq N-2$, contains the locations of
the $j$s in $\omega _{0}^{L-1}$. Each $j$-run contained in or intersected by
$\omega _{0}^{L-1}$, if any, contributes a subsequence of the same length as
the run, that is increasing ($\pi _{k_{0}+...+k_{j-1}+i},\pi
_{k_{0}+...+k_{j-1}+i}+1,\ldots $) if the run is followed by a symbol $>j$,
or \textit{decreasing} ($\pi _{k_{0}+...+k_{j-1}+i},\pi
_{k_{0}+...+k_{j-1}+i}-1,\ldots $) if the run is followed by a symbol $<j$.
This subsequences may be intertwined with other entries of the same block.

(R4) If the entries $\pi _{a}\leq L-2$ and $\pi _{b}\leq L-2$ belong to the
\textit{same block of }$\pi \in \mathcal{S}_{L}$, and $\pi _{a}$ appears on
the left of $\pi _{b}$ (i.e., $0\leq a<b\leq L-1$), then $\pi _{a}+1$
appears also on the left of $\pi _{b}+1$ (i.e., $\pi _{a}+1=\pi _{a^{\prime
}}$, $\pi _{b}+1=\pi _{b^{\prime }}$ and $0\leq a^{\prime }<b^{\prime }\leq
L-1$).

\medskip

In (R4), $\pi _{a}+1$ and $\pi _{b}+1$ may appear in the same block or in
different blocks. Let us mention at this point that (R4) implies some simple
consequences for the relative locations of increasing and decreasing
subsequences within the same block and their continuations (if any) outside
the block, but with the exception of one particular result that will be
formulated below, we will not need them in the sequel.

\medskip

\noindent \textbf{Example 6.} Take in $\{0,1,2\}^{\mathbb{N}_{0}}$ the
sequence
\begin{equation}
\left. \omega =(\right\vert _{0}\left. 2\right\vert _{1}\left. 1\right\vert
_{2}\left. 1\right\vert _{3}\left. 1\right\vert _{4}\left. 2\right\vert
_{5}\left. 2\right\vert _{6}\left. 0\right\vert _{7}\left. 0\right\vert
_{8}\left. 1\right\vert _{9}\left. 1\right\vert _{10}\left. 0\right\vert
_{11}\left. 0\right\vert _{12}\left. 2\right\vert _{13}\left. 2\right\vert
\left. 2\right\vert 1\ldots ),  \label{example}
\end{equation}%
where $\left. {}\right\vert _{k}b$ indicates that the entry $b\in \{0,1,2\}$
is at place $k$. Then $\omega $ defines the order pattern%
\begin{equation*}
\pi =[6,10,7,11;9,8,1,2,3;5,0,4,13,12]\in \mathcal{S}_{14},
\end{equation*}%
where the first block, $\pi _{0}^{3}=6,10,7,11$, is set by the $k_{0}=4$
symbols $0$ in $\omega _{0}^{13}$, which appear grouped in two runs, $\omega
_{6}^{7}$ and $\omega _{10}^{11}$ (note the two increasing subsequences $6,7$
and $10,11$ in this block); the intermediate block, $\pi _{4}^{8}=9,8,1,2,3$%
, comes from the $k_{1}=5$ symbols $1$ in $\omega _{0}^{13}$, grouped also
in two runs, $\omega _{1}^{3}$, followed by the symbol $2>1$, and $\omega
_{8}^{9}$, followed by the symbol $0<1$ (note the corresponding increasing
subsequence $1,2,3$, and decreasing subsequence $9,8$, in this block);
finally, the last block $\pi _{9}^{13}=5,0,4,13,12$ accounts for the $%
k_{2}=5 $ appearances of the symbol $2$ in $\omega _{0}^{13}$ (the
decreasing subsequences $5,4$ and $13,12$ come from the runs $\omega
_{4}^{5} $ and $\omega _{12}^{13}$, respectively, where $\omega _{12}^{13}$
is the intersection within $\omega _{0}^{13}$ of a longer $2$-run). $\square
$

Observe that two sequences $\omega $, $\omega ^{\prime }$ with $\omega
_{0}^{L-1}\neq \omega _{0}^{\prime L-1}$ may define the same order pattern
of length $L$, while two sequences $\omega $, $\omega ^{\prime }$ with $%
\omega _{0}^{L-1}=\omega _{0}^{\prime L-1}$ may define different order
patterns of length $L$ (depending on $\omega _{L-1},...,$ and $\omega
_{L-1}^{\prime },...$).

\medskip

\noindent \textbf{Proposition 3.} The one-sided shift on $N\geq 2$ symbols
has no forbidden patterns of length $L\leq N+1$.

\medskip

\noindent \textit{Proof}. First of all, note that if $\omega=(\omega_0,%
\omega_1,\omega_2,\dots)$ is of type $\pi =[\pi _{0},\pi _{1},\dots,\pi
_{N}] $, then the point $\bar\omega=(N-1-\omega_0,N-1-\omega_1,N-1-\omega_2,%
\dots)$ is of type $\pi_{mirrored}=[\pi _{N},\pi_{N-1},\dots,\pi
_{1},\pi_{0}]$.

Given $\pi =[\pi _{0},\pi _{1},\dots,\pi _{N}]$, we can therefore assume,
without loss of generality, that $\pi_0<\pi_N$. Consider two cases.

\begin{itemize}
\item If $\pi _{N}\neq N$, then there is some $l\in \{1,2,\dots ,N-1\}$ such
that $\pi _{l}=N$. In this case, the point $\omega =(\omega _{0},\omega
_{1},\dots )\in \{0,1,\dots ,N-1\}^{\mathbb{N}_{0}}$, where
\begin{eqnarray}
&&\omega _{\pi _{0}}=0,\ \omega _{\pi _{1}}=1,\ \dots ,\ \omega _{\pi
_{l-1}}=l-1,\ \omega _{\pi _{l}}=l-1,\ \omega _{\pi _{l+1}}=l,\ \dots ,
\notag \\
&&\omega _{\pi _{N-1}}=N-2,\ \omega _{\pi _{N}}=N-1,\ \omega _{N+1}=\omega
_{N+2}=N-1,  \notag
\end{eqnarray}%
is of type $\pi $. Indeed, it is enough to note that
\begin{equation*}
\Sigma ^{\pi _{l-1}}(\omega )=(l-1,\omega _{\pi _{l-1}+1},\dots
)<(l-1,N-1,N-1,\dots )=\Sigma ^{N}(\omega )=\Sigma ^{\pi _{l}}(\omega ).
\end{equation*}

\item If $\pi_N=N$, let us first assume that $\pi_0\neq0$. Then there is $%
k\in \{1,2,\dots,N-1\}$ such that $\pi_{k}+1=\pi_{0}$. In this case, the
point $\omega =(\omega _{0},\omega _{1},\dots)\in \{0,1,\dots,N-1\}^{\mathbb{%
N}_{0}}$, where
\begin{eqnarray}
&\omega _{\pi _{0}}=0,\ \omega_{\pi_1}=1,\ \dots ,\ \omega _{\pi
_{k-1}}=k-1,\ \omega _{\pi _{k}}=k,\ \omega _{\pi _{k+1}}=k,\ \omega _{\pi
_{k+2}}=k+1,\ \dots,  \notag \\
&\omega _{\pi _{N-1}}=N-2,\ \omega _{\pi _{N}}=N-1,\ \omega_{N+1}=N-1,
\notag
\end{eqnarray}%
is of type $\pi $. This is clear because
\begin{equation*}
\Sigma ^{\pi _{k}}(\omega )=(k,0,\dots)<(k,\omega _{\pi
_{k+1}+1},\dots)=\Sigma ^{\pi _{k+1}}(\omega ).
\end{equation*}

In the case that $\pi_0=0$, then there is $l\in \{1,2,\dots,N-1\}$ such that
$\pi_{l}=N-1$. Now the point $\omega =(\omega _{0},\omega _{1},\dots)\in
\{0,1,\dots,N-1\}^{\mathbb{N}_{0}}$, where
\begin{eqnarray}
&\omega_{\pi_0}=0,\ \omega_{\pi_1}=1,\ \dots,\ \omega_{\pi_{l-1}}=l-1,\
\omega_{\pi_l}=l-1,\ \omega_{\pi_{l+1}}=l,\ \dots,  \notag \\
&\omega _{\pi _{N-1}}=N-2,\ \omega_{\pi_N}=N-1,  \notag
\end{eqnarray}%
is of type $\pi $, since
\begin{equation*}
\Sigma ^{\pi _{l-1}}(\omega
)=(l-1,\omega_{\pi_{l-1}+1},\dots)<(l-1,N-1,\dots)=\Sigma ^{N-1}(\omega
)=\Sigma ^{\pi _{l}}(\omega ). \ \square
\end{equation*}
\end{itemize}

\medskip

\noindent \textbf{Proposition 4.} The shift on $N$ symbols has forbidden
patterns of length $L\geq N+2$.

\medskip

\noindent \textit{Proof}. We need only to prove the existence of forbidden
patterns of length $L=N+2$, since then their outgrowth patterns will provide
forbidden patterns of arbitrary length $L>N+2$.

Consider first the case of an \textit{even} number of symbols $%
\{0,1,...,N-1\}$, $N=2l+1$, $l\geq 1$. We claim that the `spiralling' pattern%
\begin{equation}
\pi =[2l+1,2l-1,...,3,1,0,2,...,2l,2l+2]\in \mathcal{S}_{2l+3}=\mathcal{S}%
_{N+2}  \label{pii}
\end{equation}%
is forbidden. Indeed, the central components $\pi _{l}=1$ and $\pi _{l+1}=0$
may not be in the same block, otherwise the restriction (R4) would be
violated ($2$ should be on the left of $1$). Thus we separate them with a
semicolon:%
\begin{equation*}
\pi =[2l+1,2l-1,...,3,1;0,2,...,2l,2l+2].
\end{equation*}%
Likewise, $\pi _{l+1}=0$ and $\pi _{l+2}=2$ may not be in the same block
(otherwise, according to (R4) $1$ should be on the left of $3$), hence we
separate them with a second semicolon:%
\begin{equation*}
\pi =[2l+1,2l-1,...,...,3,1;0;2,...,2l,2l+2].
\end{equation*}%
The procedure continues along alternating, outgoing directions, considering
each time pairs of consecutive components of $\pi $ in an exhaustive way: $%
\pi _{l-1}=3$ and $\pi _{l}=1$ in the 3rd step, $\pi _{l+2}=2$ and $\pi
_{l+3}=4$ in the fourth step, etc.. In the $k$th step we pick up (a) $\pi
_{l+\nu }=k-2$ and $\pi _{l+\nu +1}=k$ if $k=2\nu $, $\nu \geq 1$, or (b) $%
\pi _{l-\nu }=k$ and $\pi _{l-\nu +1}=k-2$ if $k=2\nu +1$, $\nu \geq 1$, and
conclude as before that we need to separate the corresponding pair with a $k$%
th semicolon (to put them in different blocks) in order not to violate (R4),
since $k-1$ and $k+1$ appear always in the wrong order. By the time that,
after completing the ($N-1$)th step, we arrive at the leftmost pair $\pi
_{0}=2l+1$, $\pi _{1}=2l-1$, we have already used up all the $N-1$
semicolons we have. However, this leftmost pair also violates (R4) because $%
2l+2=\pi _{N+1}$ appears on the right of $2l=\pi _{N}$. This proves that $%
\pi $ is forbidden.

Suppose now that the number of symbols is \textit{odd}: $N=2l$, $l\geq 1$.
In this case we claim that%
\begin{equation}
\tau =[2l+1,2l-1,...3,1,0,2,...,2l-2,2l]\in \mathcal{S}_{2l+2}=\mathcal{S}%
_{N+2}  \label{tau}
\end{equation}%
is forbidden. We start again considering the central components $\tau _{l}=1$
and $\tau _{l+1}=0$, to conclude that they may not be in the same block
because $2$ is on the right of $1$, violating otherwise the restriction
(R4). Thus we separate them with a semicolon:%
\begin{equation*}
\tau =[2l+1,2l-1,...,3,1;0,2,...,2l-2,2l].
\end{equation*}%
The proof continues exactly as before, except that now, after completing the
($N-1$)th step and thus having already used up $N-1$ semicolons, we arrive
at the rightmost pair $\tau _{N}=2l-2$, $\tau _{N+1}=2l$. But this pair
violates (R4) because $2l-1=\tau _{1}$ appears on the right of $2l+1=\tau
_{0}$. This proves that $\tau $ is forbidden and completes the proof. $%
\square $

\bigskip

From Proposition 3 and the proof of Proposition 4 it follows that the order
pattern (\ref{pii}) if $N$ is even, or (\ref{tau}) if $N$ is odd, is a
forbidden \textit{root} pattern of the shift on $N$ symbols. We turn next to
the question, whether there exist also forbidden root patterns of lengths $%
L>N+2$.

Consider a partition of the sequence $0,1,...,L-1$ of the form%
\begin{equation}
p_{1}<p_{2}<...<p_{d}<...<p_{D},  \label{seg1}
\end{equation}%
where \
\begin{equation}
p_{d}=e_{d},e_{d}+1,...,e_{d}+h_{d}-1,  \label{seg2}
\end{equation}%
$1\leq d\leq D$, $D\geq 2$, with (i) $h_{d}\geq 1$, $h_{1}+...+h_{D}=L$, and
(ii) $e_{d}+h_{d}=e_{d+1}$ for $1\leq d\leq D-1$, i.e., the \textit{follower}
of $p_{d}$, $e_{d}+h_{d}$, is the first element of $p_{d+1}$, $e_{d+1}$. We
call (\ref{seg1}) a partition of $0,1,...,L-1$ in $D$ segments, (\ref{seg2})
an \textit{increasing segment} and denote by $\overleftarrow{p_{d}}$ the
\textit{decreasing} or \textit{reversed segment}
\begin{equation*}
\overleftarrow{p_{n}}=e_{d}+h_{d}-1,...,e_{n}+1,e_{n}.
\end{equation*}%
We also call $e_{n}$ the first element of $\overleftarrow{p_{n}}$ and $%
e_{n+1}$ the follower of $\overleftarrow{p_{n}}$.

In the proof of the existence of forbidden root patterns below (Lema 1 and
Proposition 5) we are going to use the following straightforward consequence
of restriction (R4) (that we will hence also refer to as (R4)):\textit{\ The
follower }$($\textit{if any}$)$\textit{\ of an increasing segment }$p_{n}$%
\textit{\ }$($\textit{correspondingly, decreasing segment }$\overleftarrow{%
p_{n}})$\textit{\ in an allowed pattern }$\pi $\textit{\ appears always to
the right of }$p_{n}$\textit{\ }$($\textit{correspondingly, to the left of }$%
\overleftarrow{p_{n}})$.

\medskip

\noindent \textbf{Definition.} Given a partition (\ref{seg1}) of $%
0,1,...,L-1 $ in segments, we call an order pattern of the form%
\begin{equation}
\pi =[...,\overleftarrow{p_{3}},\overleftarrow{p_{1}},p_{2},p_{4},...],
\label{pinormal}
\end{equation}%
or its \textit{mirrored pattern}%
\begin{equation}
\pi _{mirrored}=[...,\overleftarrow{p_{4}},\overleftarrow{p_{2}}%
,p_{1},p_{3},...],  \label{pimirror}
\end{equation}%
a \textit{spiralling pattern} of length $L$.

\medskip

Observe that the relation between partitions of $0,1,...,L-1$ in segments
and spiralling patterns of length $L$ is one-to-one except when $p_{1}=0$ ($%
h_{1}=1$). In this case, $\overleftarrow{p_{1}},p_{2}=0,1,...,e_{2}+h_{2}-1$
can be taken for $p_{1}^{\prime }\equiv 0,1,...,e_{2}+h_{2}-1$ ($%
h_{1}^{\prime }=h_{2}+1$).

\medskip

\noindent \textbf{Lemma 1.} If $N\geq 2$ is the number of symbols and $\pi $
is a spiralling pattern with $D$ segments and $h_{1}\geq 2$ (i.e., $%
p_{1}=0,1,...$), then $\pi $ is forbidden if (a) $D\geq N$ and $h_{D}\geq 2$
or (b) $D\geq N+1$ and $h_{D}=1$; otherwise it is allowed.

\medskip

The first part of this proposition generalizes Proposition 4. Indeed, the
order patterns (\ref{pii}) and (\ref{tau}) correspond to the `minimal'
spiralling pattern in case (b): $p_{1}=0,1$ and $p_{d}=d$ for $2\leq d\leq
N+1$.

\medskip

\noindent \textit{Proof}. Consider the spiralling pattern (\ref{pinormal}).
The proof that such $\pi $ is forbidden proceeds formally as in Proposition
4, starting again with the central segment $\overleftarrow{p_{1}}%
=e_{1}+h_{1}-1,...,1,0$ (first semicolon). From here on, three possibilities
can occur that we illustrate in a general step of even order. (i) If $%
p_{2\nu }$ consists of more than one element (i.e., $h_{2\nu }\geq 2$), then
we apply (R4) to $p_{2\nu }$ to conclude that we need a semicolon between $%
e_{2\nu }+h_{2\nu }-2$ and $e_{2\nu }+h_{2}-1$ (since the follower of $%
p_{2\nu }$, i.e., the first entry of $\overleftarrow{p_{2\nu +1}}$, is on
the wrong side). (ii) If $p_{2\nu }$ consists of one element ($h_{2\nu }=1$)
and $p_{2\nu -2}$ consists of more than one element ($h_{2\nu -2}\geq 2$),
then we apply (R4) to the pair $p_{2\nu }=e_{2\nu }$ and $e_{2\nu
-2}+h_{2\nu -2}-1$, the last element of $p_{2\nu -2}$, which has been
separated with a semicolon from the rest of elements in $p_{2\nu -2}$ two
steps earlier. (iii) If both $p_{2\nu }$ and $p_{2\nu -2}$ consist of a
single element ($h_{2\nu }=h_{2\nu -2}=1$), apply (R4) to the pair $p_{2\nu
-2}=e_{2\nu -2}<p_{2\nu }=e_{2\nu }$ to infer the need for a semicolon
separating them (since $e_{2\nu -2}+1=e_{2\nu -1}$, the first element of $%
\overleftarrow{p_{2\nu -1}}$, is on the right of $e_{2\nu }+1=e_{2\nu +1},$
the first element of $\overleftarrow{p_{2\nu +1}}$). As a general rule, we
need one semicolon per segment $p_{2\nu }$ or $\overleftarrow{p_{2\nu +1}}$,
as long as there are still a posterior segment $\overleftarrow{p_{2\nu +1}}$
or $p_{2\nu +2}$, respectively, on the `wrong' side.

Following in this way, we run out of semicolons ($N-1$ at most) after having
considered the segment $p_{N-1}$. If $D\geq N$ and $h_{N}\geq 2$, then $%
p_{N} $ will violate (R1) if $N$ is odd or (R2) if $N$ is even. If $h_{N}=1$
but $D\geq N+1$, then the segment $p_{N+1}$ will be on the wrong side of $%
p_{N}$ and the pattern will not comply with (R4).

The proof for $\pi _{mirrored}$, Eq. (\ref{pimirror}), is completely
analogue.

Also, the procedure above shows how to decompose any spiralling pattern into
well-formed (i.e., complying with (R1)-(R4)) blocks. The central block (in
the case (\ref{pinormal})) is of the form $e_{1}+h_{1}-2,...,1,0$ if $%
h_{2}=1 $, or $e_{1}+h_{1}-2,...,1,0,2,..,e_{2}+h_{2}-2$ if $h_{2}\geq 2$.
Each block on the right side of the central block is of the form $e_{2\nu
}+h_{2\nu }-1$ if $h_{2\nu +2}=1$ or $2\nu =D$ (rightmost block), or it has
the form $e_{2\nu }+h_{2\nu }-1,e_{2\nu +2},...,e_{2\nu +2}+h_{2\nu +2}-2$
if $h_{2\nu +2}\geq 2$. Each block on the left side of the central block is
of the form $e_{2\nu -1}+h_{2\nu -1}-1$ if $h_{2\nu +1}=1$ or $2\nu -1=D$
(leftmost block), or it has the form $e_{2\nu +1}+h_{2\nu +1}-2,...,e_{2\nu
+1},e_{2\nu -1}+h_{2\nu -1}-1$ if $h_{2\nu +1}\geq 2$. If $N$, the number of
symbols, is equal to or greater than the number of resulting blocks ($D+1$
if $h_{D}\geq 2$, and $D$ if $h_{D}=1$), one can readily write down
sequences $\omega \in \{0,1,...,N-1\}^{\mathbb{N}_{0}}$ of type $\pi $. $%
\square $

\medskip

\noindent \textbf{Example 7.} As illustration of the procedure used in the
proof of Proposition 5, consider the spiralling pattern%
\begin{equation*}
\pi =[9,8,7,5,2,1,0,3,4,6,10,11]\in \mathcal{S}_{12}.
\end{equation*}%
Here $p_{1}=0,1,2$, $p_{2}=3,4$, $p_{3}=5$, $p_{4}=6$, $p_{5}=7,8,9$ and $%
p_{6}=10,11$. The following scheme summarizes the steps of the decomposition
of $\pi $ into well-formed blocks:%
\begin{equation*}
\begin{array}{lllllll}
\overleftarrow{p_{1}}=2,1,0 & \rightarrow  & 2;1,0 &  & p_{2}=3,4 &
\rightarrow  & 3;4 \\
\overleftarrow{p_{3}},2=5,2 & \rightarrow  & 5;2 &  & 4,p_{4}=4,6 &
\rightarrow  & 4;6 \\
\overleftarrow{p_{5}}=9,8,7 & \rightarrow  & 9;8,7 &  & p_{6}=10,11 &
\rightarrow  & 10;11%
\end{array}%
\end{equation*}%
Hence,%
\begin{equation*}
\pi =[9;8,7,5;2;1,0,3;4;6,10;11]
\end{equation*}%
Since the decomposition consists of $7$ blocks, $\pi $ is allowed if $N\geq 7
$, in compliance with Lemma 1 (a) with $D=6$ and $h_{6}=2$. For instance,
any sequence $\omega \in \{0,1,...,6\}^{\mathbb{N}_{0}}$ such that%
\begin{equation*}
\omega _{0}^{11}\equiv \omega _{0},...,\omega _{11}=3,3,2,3,4,1,5,1,1,0,5,6
\end{equation*}%
is of type $\pi $.

\medskip

\noindent \textbf{Proposition 5.} For every $L\geq N+2$, the one-sided shift
on $N$ symbols has forbidden root patterns of length $L$.

\medskip

\noindent \textit{Proof}. If $L=N+2$, we know already (Proposition 3 and 4)
that the spiralling pattern (\ref{pii}) if $N$ is even, or (\ref{tau}) if $N$
is odd, is a forbidden root pattern. Thus, assume $L>N+2$ and consider the
following partition of $0,1,...,L-1$ in $N$ segments:%
\begin{equation*}
0,1<p_{2}<...<p_{N-1}<L-2,L-1.
\end{equation*}%
(i.e., $h_{1}=h_{N}=2$). We claim that the spiralling pattern%
\begin{equation}
\pi ^{\ast }=[\overleftarrow{p_{N-1}},...,\overleftarrow{p_{3}}%
,1,0,p_{2},...,p_{N-2},L-2,L-1]  \label{roote}
\end{equation}%
if $N$ is even, or%
\begin{equation}
\tau ^{\ast }=[L-1,L-2,\overleftarrow{p_{N-2}}...,\overleftarrow{p_{3}}%
,1,0,p_{2},...,p_{N-1}],  \label{tauasterisk}
\end{equation}%
if $N$ is odd, and their corresponding mirrored patterns, are forbidden root
patterns. Only the first case will be analyzed here, the proof being
completely analogue in the second case and for the mirrored patterns.

That (\ref{roote}) is forbidden follows readily from Lemma 1 (a). To prove
next that $\pi ^{\ast }$ is a forbidden root pattern, we need to show that
it is not the outgrowth pattern of any forbidden pattern of shorter length.
Remember that given a forbidden pattern%
\begin{equation*}
\lbrack \pi _{0},...,\pi _{L-2}]\in \mathcal{S}_{L-1},
\end{equation*}%
its outgrowth patterns of length $L$ have the form (\textit{Group A})%
\begin{equation*}
\lbrack L-1,\pi _{0},...,\pi _{L-2}],[\pi _{0},L-1,...,\pi _{L-2}],...,[\pi
_{0},...,\pi _{L-2},L-1],
\end{equation*}%
or the form (\textit{Group B})%
\begin{equation*}
\lbrack 0,\pi _{0}+1,...,\pi _{L-2}+1],[\pi _{0}+1,0,...,\pi
_{L-2}+1],...,[\pi _{0}+1,...,\pi _{L-2}+1,0].
\end{equation*}

There are two possibilities. Suppose first that $\pi ^{\ast }$ is an
outgrowth forbidden pattern of Group A. Then deleting the entry $L-1$ yields
the spiralling pattern%
\begin{equation*}
\lbrack \overleftarrow{p_{N-1}},...,\overleftarrow{p_{3}}%
,1,0,p_{2},...,p_{N-2},L-2],
\end{equation*}%
which is allowed on account of having $N$ segments, $h_{1}=2$, and a last
segment of length $1$ (Lemma 1 (b)).

Thus, suppose that $\pi ^{\ast }$ is an outgrowth forbidden pattern of Group
B. Then deleting the entry $0$ and subtracting 1 from the remaining entries,
we get the pattern%
\begin{equation}
\lbrack \overleftarrow{p_{N-1}^{\prime }},...,\overleftarrow{p_{3}^{\prime }}%
,0,p_{2}^{\prime },...,p_{N-2}^{\prime },L-3,L-2],  \label{B}
\end{equation}%
where $p_{d}^{\prime }=$ $e_{d}-1,...,e_{d}+h_{d}-2$, $1\leq d\leq N+1$.
Since $p_{1}^{\prime }=0$ ($h_{1}^{\prime }=h_{1}-1=1$) and $p_{2}^{\prime
}=1,...$ ($h_{2}^{\prime }=h_{2}\geq 1$), we can merge $p_{1}^{\prime }$ and
$p_{2}^{\prime }$ into the new segment $p_{1}^{\prime \prime }\equiv 0,1,...$%
, so that (\ref{B}) is a spiralling pattern with $h_{1}^{\prime \prime }\geq
2$ and the following $N-1$ segments: $p_{1}^{\prime \prime },p_{3}^{\prime
},...,p_{N-1}^{\prime }$ and $p_{N}^{\prime }=L-3,L-2$. According to Lemma 1
(a), the order pattern (\ref{B}) is allowed. $\square $

\bigskip

\noindent\textbf{Remark. }Spiralling patterns of the particular form (\ref%
{roote}) or (\ref{tauasterisk}) (and the corresponding mirrored patterns)
are not, of course, the only forbidden root patterns for the shift on $N$
symbols. For instance, it can be easily checked that all patterns of length $%
L\geq 2N$ the form%
\begin{equation*}
\lbrack 1;0,3;2,5;4...,2N-3;2N-4,L-2,L-3,...,2N-2,L-1]\in \mathcal{S}_{L}
\end{equation*}%
and their mirrored patterns, are forbidden root patterns as well.

\medskip

\noindent \textbf{Corollary 1.} For every $K\geq 2$ there are maps on $[0,1]$
without forbidden patterns of length $L\leq K$.

\medskip

\noindent \textit{Proof}. Let $S_{N}=\psi \circ \Sigma \circ \psi
^{-1}:[0,1]\rightarrow \lbrack 0,1]$ be the map (\ref{sawtooth}). Since $%
\psi $ is an order-isomorphy, $S_{N}$ and $\Sigma $, the shift on $N$
symbols, have the same forbidden patterns. Therefore, if $N+1\leq K$, then $%
S_{N}$ has no forbidden patterns of length $L\leq K$ because of Proposition
3. $\square $

\medskip

It follows that \textit{there are interval maps on} $\mathbb{R}^{n}$ \textit{%
without forbidden patterns}. For example, one can decompose $[0,1]$ in
infinite many half-open intervals (of vanishing length), $[0,1]=\cup
_{N=2}^{\infty }I_{N}$ and define on each $I_{N}$ a properly scaled version
of $S_{N}$, $\tilde{S}_{N}:I_{N}\rightarrow I_{N}$. In $\mathbb{R}^{2}$ one
can perform the said decomposition along the $1$-axis and define on $%
I_{N}\times \lbrack 0,1]$ the function $(\tilde{S}_{N}$, Id$).$ Now, Eq. (%
\ref{TopEntropy2}) shows that adding some natural assumption, like piecewise
monotony, can make all the difference.

\medskip

\section{\noindent \textbf{Order patterns and two-sided shifts}}

Consider now the bisequence space, $\{0,1,$ $...,N-1\}^{\mathbb{Z}}$,
endowed with the lexicographical (or product) order. With the notation $%
\omega _{-}$ for the \textit{left sequence} $(\omega _{-n})_{n\in \mathbb{N}%
} $ of $\omega \in \{0,1,$ $...,N-1\}^{\mathbb{Z}}$ and $\omega _{+}$ for
the \textit{right sequence} $(\omega _{n})_{n\in \mathbb{N}_{0}}$, we have
\begin{equation*}
\omega <\omega ^{\prime }\;\Leftrightarrow \left\{
\begin{array}{l}
\;\omega _{+}<\omega _{+}^{\prime } \\
\text{\mbox{or}} \\
\omega _{-}<\omega _{-}^{\prime }\text{ \mbox{if} }\omega _{+}=\omega
_{+}^{\prime }%
\end{array}%
\right. \text{ ,}
\end{equation*}%
where $<$ between right (resp. left) sequences denotes lexicographical order
in $\{0,1,$ $...,N-1\}^{\mathbb{N}_{0}}$ (resp. $\{0,1,$ $...,N-1\}^{\mathbb{%
N}}$). Thus, the lexicographical order for bisequences is defined most of
the time by the right sequences of the points being compared, except when
they coincide, in which case the order is defined by their left sequences.
If we map $\{0,1,$ $...,N-1\}^{\mathbb{Z}}$ onto $[0,1]\times \lbrack
0,1]\equiv \lbrack 0,1]^{2}$ via
\begin{equation*}
(\omega _{-},\omega _{+})\mapsto \left( \sum_{n=1}^{\infty }\omega
_{-n}N^{-n},\sum_{n=0}^{\infty }\omega _{n}N^{-(n+1)}\right) ,
\end{equation*}%
we find that lexicographical order in $\{0,1,$ $...,N-1\}^{\mathbb{Z}}$
corresponds to lexicographical order in $[0,1]^{2}$, which results thereby
foliated into a continuum of copies of $([0,1],<)$. In order for this map to
be one-to-one, we have to exclude the countable set $\overline{\overline{%
\mathcal{N}}}$ of all bisequences terminating in an infinite string of $%
(N-1) $s in either direction.

In relation with the order patterns defined by the orbits of two-sided
sequences,

\begin{equation*}
\Sigma ^{i}(\omega )<\Sigma ^{j}(\omega )\;\Leftrightarrow \left\{
\begin{array}{l}
\;(\omega _{i},\omega _{i+1},...)<(\omega _{j},\omega _{j+1},...) \\
\text{\mbox{or}} \\
(\omega _{i-1},\omega _{i-2},...)<(\omega _{j-1},\omega _{j-2},...)\text{ %
\mbox{if} }(\omega _{i},\omega _{i+1},...)=(\omega _{j},\omega _{j+1},...)%
\end{array}%
\right. \text{ ,}
\end{equation*}%
where $i,j\geq 0$, $i\neq j$. It follows that the `exceptional' condition $%
(\omega _{i},\omega _{i+1},...)=(\omega _{j},\omega _{j+1},...)$ occurs if
and only if $\Sigma ^{\left\vert i-j\right\vert }(\omega _{+})=\omega _{+}$,
i.e., when the right sequence $\omega _{+}$ of $\omega \in \{0,1,$ $%
...,N-1\}^{\mathbb{Z}}$ is periodic from the entry $\min \{i,j\}$ on with
period $p=\left\vert i-j\right\vert $.

\medskip

\noindent \textbf{Proposition 6.} The two-sided shift on $N\geq 2$ symbols
has no forbidden patterns of length $L\leq N-1$ and has forbidden patterns
for $L\geq N+2$.

\medskip

\noindent \textit{Proof}. The one-sided sequence $\omega _{+}\in \{0,1,$ $%
...,N-1\}^{\mathbb{N}_{0}}$ defines an order pattern $\pi $ of length $L$,%
\begin{equation*}
\Sigma ^{\pi _{0}}(\omega _{+})<\Sigma ^{\pi _{1}}(\omega _{+})<...<\Sigma
^{\pi _{L-1}}(\omega _{+}),
\end{equation*}%
if and only if the two-sided sequences $\omega =(\omega _{-},\omega _{+})$,
with $\omega _{-}$ $\in \{0,1,$ $...,N-1\}^{\mathbb{N}}$ arbitrary, define
the same order pattern. $\square $

\medskip

\noindent \textbf{Example 8.} Let $I^{2}=[0,1]\times \lbrack 0,1]$ endowed
with the induced Lebesgue measure $\lambda $ and $B:I^{2}\rightarrow I^{2}$
the $\lambda $-invariant \textit{baker's map},%
\begin{equation*}
B(x,y)=\left\{
\begin{array}{lc}
(2x,\frac{1}{2}y), & 0\leq x<\frac{1}{2}, \\
(2x-1,\frac{1}{2}y+\frac{1}{2}), & \frac{1}{2}\leq x\leq 1.%
\end{array}%
\right.
\end{equation*}%
A generating partition of $(I^{2},\lambda ,B)$ is: $A_{0}=[0,\frac{1}{2}%
)\times \lbrack 0,1]$ and $A_{1}=[\frac{1}{2},1]\times \lbrack 0,1]$. For $%
\Sigma $ take the two-sided $(\frac{1}{2},\frac{1}{2})$-Bernoulli shift.
Then $B$ and $\Sigma $ are isomorphic via the $\lambda $-invariant coding
map $\Phi :I^{2}\rightarrow \{0,2\}^{\mathbb{Z}}\backslash $ $\overline{%
\overline{\mathcal{N}}}$, given by%
\begin{equation*}
\Phi (x)=(...,\omega _{-1},\omega _{0},\omega _{1},...),
\end{equation*}%
where $\omega _{n}=a_{n}$ if $B^{n}(x)\in A_{a_{n}}$, $n\in \mathbb{Z}$.
Since $\Phi $ preserves order (in fact, $\Phi $ is the inverse of the
order-preserving map $(\omega _{-},\omega _{+})\mapsto (\sum_{n=0}^{\infty
}\omega _{-n}2^{-(n+1)},\sum_{n=1}^{\infty }\omega _{n}2^{-n})$, sequences
ending with $\bar{1}$ excluded), we conclude that the baker's transformation
has no forbidden patterns of length $\leq 3$. $\square $

\medskip

\noindent {\small ACKNOWLEDGEMENTS. We are very thankful to the referees for
their valuable comments. J.M.A. was financially supported by the Spanish
Ministry of Education and Science (Grant MTM2005-04948).}

\end{document}